\documentclass[11pt,a4paper]{amsart}

\usepackage{ae} 
\usepackage[cp1250]{inputenc}
\usepackage{amsmath}
\usepackage{amssymb}
\usepackage{color}
\usepackage[colorlinks]{hyperref}
\usepackage[slovene, english]{babel}
\usepackage[dvips]{graphicx}
\usepackage[format=plain,labelfont=bf,up,textfont=it,up]{caption}

\def\notation{{\bf Notation. }}
\def\proof{{\bf Proof. }}

\def\remark{{\hskip-6mm \bf Remark. }}

\newtheorem{dfn}{Definition}
\newtheorem{prob}[dfn]{Problem}
\newtheorem{prop}[dfn]{Proposition}
\newtheorem{thm}[dfn]{Theorem}
\newtheorem{lemma}[dfn]{Lemma}

\newtheorem{cor}[dfn]{Corollary}

\renewcommand{\thedfn}%
{\arabic{dfn}}

\renewcommand{\rightsquigarrow}%
{ \hspace{4mm}\rput{0}(0,.1){\diagup \hskip-1.3mm \searrow}\hspace{4mm} }

\newcommand{\Int}{\mathop%
{\emph{Int}}}
\newcommand{\rel}%
{ \emph{\hspace{0.1 cm}rel } }

\def\eps{\varepsilon}

\def\NN{\mathbb{N}}
\def\QQ{\mathbb{Q}}
\def\RR{\mathbb{R}}
\def\ZZ{\mathbb{Z}}

\newcommand{\rp}%
{ is relatively prime to }
\newcommand{\nrp}%
{ is not relatively prime to }

\def\di{\partial}                

\def\di{\partial}
\def\fr{\frac}

\def\a{\alpha}
\def\b{\beta}
\def\g{\gamma}

\def\f{\varphi}
\def\F{\Phi}

\renewcommand{\theenumi}%
{\roman{enumi}}
\renewcommand{\labelenumi}%
{(\theenumi)}

\input xy
\xyoption{all}

\title{Realizations of  Countable Groups as Fundamental Groups of Compacta}
\author{ \v Ziga Virk }
\thanks{The author would like to thank Jerzy Dydak and Ale\v s Vavpeti\v c for their valuable comments.}
\address{University of Tennessee, Knoxville, TN 37996}
\email{zigavirk@gmail.com}

\begin{document}

\begin{abstract}
It is an open question (Pawlikowski, \cite{Paw}) whether every finitely generated  group can be realized as a fundamental group of a compact metric space. In this paper we prove that any countable group can be realized as the fundamental group of a compact subspace of $\RR^4.$ According to theorems of Shelah \cite{SH1} (see also Pawlikowski \cite{Paw})  such space can not be locally path connected if the group is not finitely generated. This constructions complements realization of groups in the context of compact Hausdorff spaces, that was studied by Keesling and Rudyak\cite{Kee},  and  Prze\'zdziecki \cite{Prze}.
\end{abstract}

\maketitle

\section{Introduction}

It is of interest to investigate analogs of constructions and results
of algebraic topology of $CW$ complexes in the category of compact spaces. For example, S.Ferry \cite{Ferry} proved that every space
homotopy dominated by a compactum is homotopy equivalent to a compactum. Analogous result holds for $CW$ complexes but fails for finite $CW$ complexes.

In this paper we deal with the problem of creating a compact
space $X_G$ with a given fundamental group $G$.
That problem was discussed in papers \cite{SH1}, \cite{Paw},
and \cite{Kee}. Shelah \cite{SH1} proved $G$ must be finitely
generated if $G$ is countable and $X_G$ is a Peano continuum.
An alternative proof of that result was presented by Pawlikowski
\cite{Paw} who posed the reverse question:

\begin{prob}[Pawlikowski] Given a finitely generated group
$G$ is there a continuum $X_G$ such that
$\pi_1(X_G)=G$.
\end{prob}

Keesling and Rudyak \cite{Kee} addressed the case of
groups $G$ for which $X_G$ can be chosen as compact Hausdorff.
Namely, every group of non-measurable cardinality is the fundamental group of a compact space.
However, their construction yields non-metrizable and non-path connected spaces, so they posed the following question
in the electronic version of their paper:

\begin{prob}[Keesling and Rudyak]\label{KRQ} For which groups
$G$ is there a path-connected compact Hausdorff $X_G$ such that
$\pi_1(X_G)=G$?
\end{prob}
 Adam Prze\'zdziecki \cite{Prze}
answered \ref{KRQ} in affirmative for any $G$ of non-measurable cardinality. Also, he announced an example of an abelian
group of measurable cardinality that is not the fundamental group of any compact
space.

As is well-known every group $G$ admits a $CW$ complex $K_G$
satisfying $\pi_1(K_G)=G$. $K_G$ may be finite if and only if $G$ is finitely presented.

A natural idea when constructing a space with prescribed fundamental group $G$ is to realize it as a two dimensional $CW$ complex $K_G:$ $1-$cells correspond to generators of $G$ and $2-$cells correspond to relations of $G$. In  case of  general countable or even finitely generated group $K_G$ may not be metric or compact. In the case of $G$ being countable we plan to construct a compact metric space $X_G$ with the fundamental group isomorphic to $G$. The idea is to replace $1-$cells (using a variation of smallness property \cite{ZV}) by suitable spaces which will enable our construction to take place in $\RR^4$. Such replacement will allow us to make our space compact but we will lose local path connectedness and  universal property for extending maps that $K_G$ has: any homomorphism $G\to \pi_1(Y)$ induces $K_G\to Y$ for any space $Y$.

\section{Construction}

Let $G$ be any countable group with presentation as $\langle g_1, g_2, \ldots \mid r_1, r_2, \ldots \rangle.$ We will construct compact metric space $X_G$ that has $G$ as  fundamental group.

\subsection{Harmonic Vase}

The basic  step in our construction is the Harmonic Vase. It replaces $1-$cells in the construction of $K_G$.
\begin{dfn}
\label{HV}
The \textbf{Harmonic Vase} with parameters $m,p\in \RR^+$ [notation: $HV(m,p)$] is the subset of $\RR^3$ defined as the union of two sets:
\begin{itemize}
  \item the \textbf{pedestal} $B(3,0)\cap (\RR^2 \times \{0\})=\{(x,y,0)\in \RR^3, x^2+y^2\leq 9\}=\{r\leq 3, z=0\},$ and
  \item the \textbf{wall} $W(m,p)$, parameterized as
  $$
  z\in (0,m],\quad \f\in [-\pi,\pi], \quad r:= \frac{|\f|}{\pi}\sin \frac{\pi p}{z}+2
  $$
  where $(r,\f)$ are polar coordinates in $\RR^2\times \{0\}\subset \RR^3$ and $z$ is the  coordinate of $\{0\}^2\times \RR$ so that $(r,\f,z)$ are cylindric coordinates in $\RR^3$.
\end{itemize}
\end{dfn}

\begin{figure}
\includegraphics[bb=90 520 450 760]{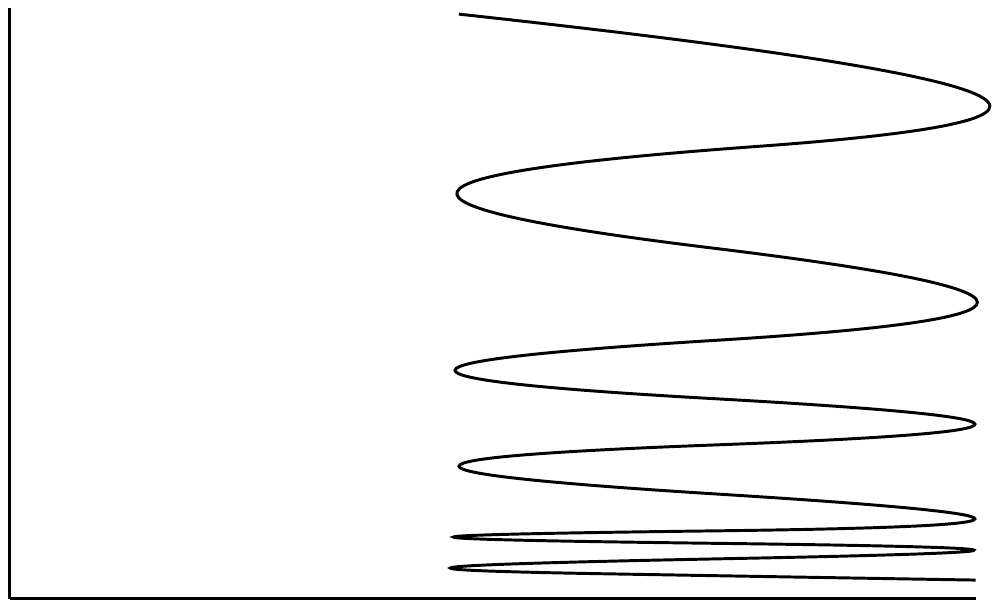}
\caption{Intersection of an $HV$ with $\f \in \{0,\pi\}$.}
\label{HVcut1}
\end{figure}

\begin{figure}
\includegraphics[bb=90 580 450 770]{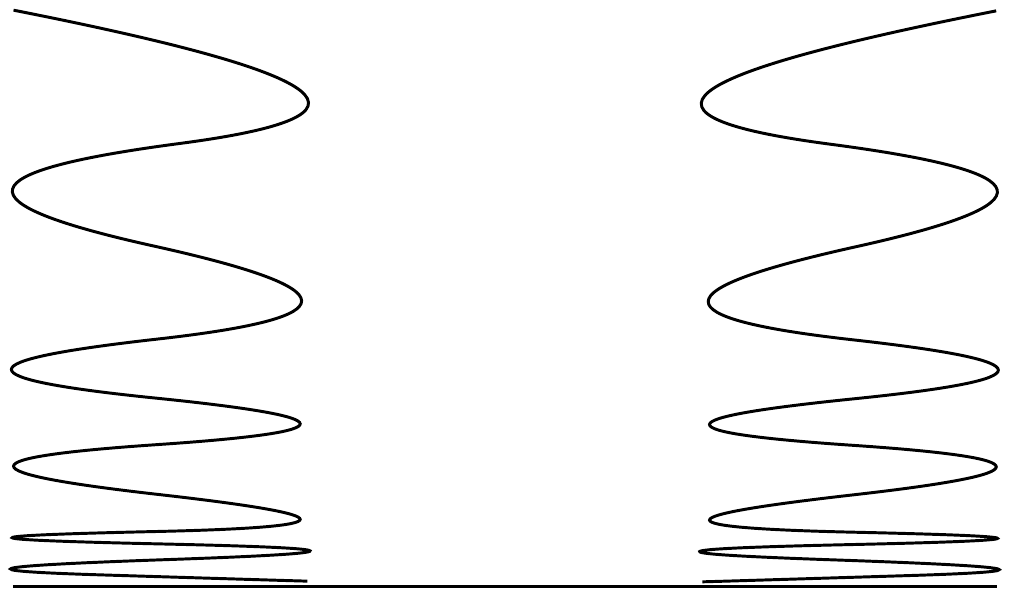}
\caption{Intersection of an $HV$ with $\f=\pm \pi/2$.}
\end{figure}

\begin{figure}
\includegraphics[scale=0.85,bb=150 210 450 610]{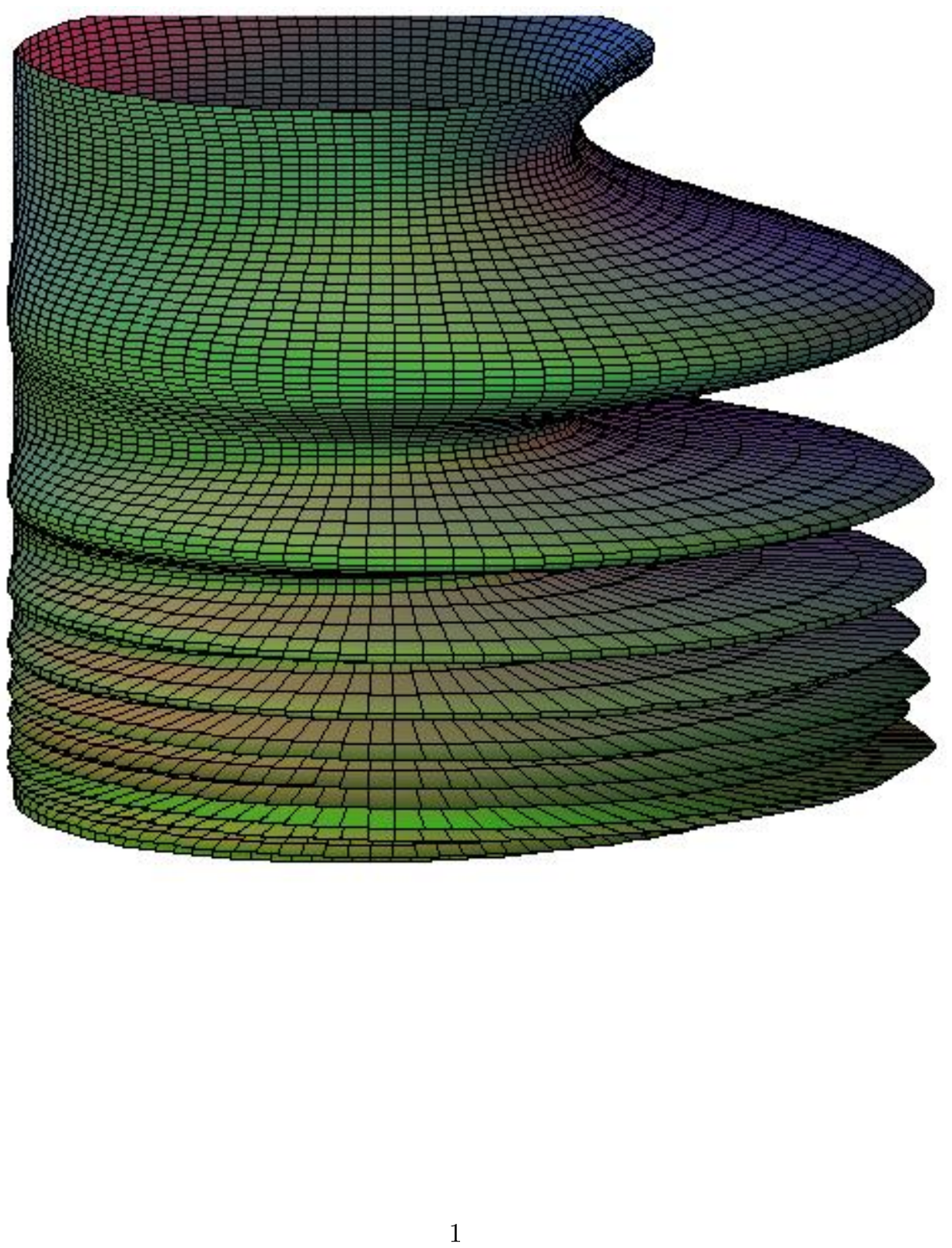}
\caption{The wall of an $HV$ drawn in Maple.}
\end{figure}

For the sake of simplicity we will use the notation $HV$ instead of $HV(m,p)$ if the parameters don't play any crucial role.

Using cylindrical coordinates, we will always assume $\f\in [-\pi,\pi]$.
To visualize the wall of $HV$ note that for any fixed $a\neq 0, a\in [-\pi,\pi]$ the intersection of the wall and halfplane $\f=a$ in $\RR^3$ is a reparameterized $\sin\frac{1}{x}$ curve. In case $a=0$  we  get a semi open line segment.

The parameter $m$ in $HV(m,p)$ is the height of the Harmonic Vase, the parameter $p$ determines a parametrization of $\sin\frac{1}{x}$ curving of the wall. We will vary both of these in our construction.

\begin{prop}\label{CompactHV}
Every $HV$ is compact.
\end{prop}

\proof
Take any Cauchy sequence $S$ in $HV$. If $S$ converges to a point in $\RR^3$ with $z=0$ then that limit point is contained in the pedestal of $HV$. If $S$ converges to a point in $\RR^3$ with $z>0$ then we may assume that all points of $S$ have $z-$coordinate at least $\eps$ for some fixed $\eps>0$. Because $HV\cap \{z\geq \eps\}$ is compact (by the definition it is the image of $[\eps,m]\times [-\pi,\pi]$ under a continuous function) and contains $S$, it also contains the limit point.
\hfill $\blacksquare$
\bigskip

When constructing the realization space $X_G$, we will use $HV$'s with various parameters. In order to combine them efficiently we have to introduce the notion of inner-curves of $HV(m,p)$.

\begin{dfn} The \textbf{inner-heights} of $HV(m,p)$ are numbers $\{h\in \RR^+\mid \sin(\frac{\pi p}{h})=-1\}$. The \textbf{inner-curves} are simple closed curves $S(m,p,c):=HV(m,p)\cap \{z=c\}$ where $c\in (0,m]$ is any height.
\end{dfn}

For any fixed $a\in (0,m]$ the orthogonal projection of the curve $S(m,p,a)$ to $\RR^2\times \{0\}$ is parameterized as
$$
R= \frac{|\F|}{\pi}\sin \frac{\pi p}{a}+2, \F\in[-\pi,\pi],
$$
where $(R,\F)$ are polar coordinates in $\RR^2$. For any two choices of $a$, these curves are smooth topological circles that  either have the only common point at $\F=0$ (in which case one of the curves lies inside the set bounded by the other curve)  or they are same hence we can talk about some of these curves  being inner or outer. The meaning of inner-heights is to provide the set of heights, for which these projections are inner-most curves.

As mentioned above, $HV$'s  replace $S^1$ in the construction of $CW$ countable group realization. It is hence important to know it's fundamental group. For this purpose we recall the definition of the universal Peano space introduced in \cite{D1}.

\begin{dfn} Let $X$ be a path connected space. The \textbf{universal Peano space} of $X$ [notation: $P(X)$] is the set $X$ equipped with a new topology, generated by all path-components of all open subsets of the existing topology on $X$. The \textbf{universal Peano map} is the natural bijection $P(X)\to X$.
\end{dfn}

The name "universal Peano map" refers to the universal map lifting property for locally path connected spaces:

\begin{prop} Let $Y$ be a locally path connected space. Then any map $Y\to X$ uniquely lifts to a map $Y\to P(X)$.
\end{prop}
$$
\xymatrix{&P(X) \ar[d]\\
Y \ar@{-->}[ur] \ar[r]&X
}
$$
\bigskip

The proof of this proposition is easy and provided in \cite{D1}. Note that if $Y$ is locally path connected then so is $Y\times I$ (where $I:=[0,1]$) which yields  the following corollary.

\begin{cor} \label{MapsToPeano} Let $Y$ be   locally path connected space and let $X$ be path connected space.
\begin{enumerate}
  \item The set of homotopy classes of maps $[Y,X]$ is in a natural bijection with $[Y,P(X)].$
  \item  The set of homotopy classes of maps $[Y,X]_\bullet$ in the pointed category is in a natural bijection with $[Y,P(X)]_\bullet.$
  \item $\pi_k(X)=\pi_k(P(X))$, for all $k\in \ZZ^+$.
\end{enumerate}
\end{cor}

\begin{prop} \label{PeanifHV}
For every choice of parameters $m,p$ one has $\pi_*(HV(m,p))=\pi_*(S^1)$. Moreover, the inclusion of the top edge of $HV$ into $HV$ is a weak homotopy equivalence.
\end{prop}

\proof To simplify the notation in this proof we will use notation $HV$ instead of $HV(m,p)$.

The crucial step is to extract space $P(HV)$. We have to consider four different types of points.

\begin{enumerate}
  \item Every point of  the  wall of $HV$ has arbitrarily small simply connected neighborhood as  the wall itself is homeomorphic to $S^1\times (0,1]$. Hence  the topology of $P(HV)$ at those points is no different from the topology in $HV$.

  \item The point $(z=0,\f=0,r=2)$ (we will mark that point with $x_0$) also has arbitrarily small simply connected neighborhood, which is a bit harder to see. Let $\eps>0$ be sufficiently small and consider neighborhood $U_\eps$ of $x_0$ in $HV$ that contains all points with $z<\eps, |\f|<\eps, 2+\eps >r>2-\eps.$ We will show that any point $x\in U_\eps$ can be connected to $x_0$ by a path which is enough for  the proof of our claim.

      If  the $z-$coordinate of $x$ equals $0$ then $x$ can be connected to $x_0$ because the pedestal (that contains both $x_0$ and $x$) is locally path connected.

      If  the $z-$coordinate of $x$ equals $h\neq 0$ then $U_\eps \cap \{z=h\}$ is an open arc, containing path from $x$ to a point $(z=h,\f=0,r=2)$. This point can be connected to $x_0$ in $H_\eps$ by a straight line segment.

  \item The points on pedestal with $r>\frac{|\f|}{\pi}+2$ or $r<-\frac{|\f|}{\pi}+2$ are not limit points of the wall. Hence they all have arbitrarily small simply connected neighborhood.

  \item The points on pedestal with $\frac{|\f|}{\pi}+2 \geq r \geq-\frac{|\f|}{\pi}+2$ are limit points of the wall. But any point $x$, other that $x_0$ has  a neighborhood $U_x$  small enough, so that  the  path component of $U_x$ containing $x$ lies in  the  pedestal.
\end{enumerate}

Summing up, the only change of topology happens at the points in $(iv)$ which become separated from the wall: they have neighborhoods contained only in the pedestal. This means that $P(HV)$  is homeomorphic to a wedge $B^2 \vee \Big((S^1\times (0,1]) \cup \{x_0\}\Big)$ where $B^2$ corresponds ti  pedestal, wedge point represents $x_0$ and $S^1\times (0,1]$ corresponds to the wall of $HV$.  There is a strong deformation retraction $P(HV)\to S^1$. Using  homotopic equivalence $P(HV)\simeq S^1$ and \ref{MapsToPeano}we get $\pi_*(HV)\cong \pi_*(S^1)$.
\hfill $\blacksquare$
\bigskip

\begin{figure}
    \includegraphics[bb=100 590 450 780]{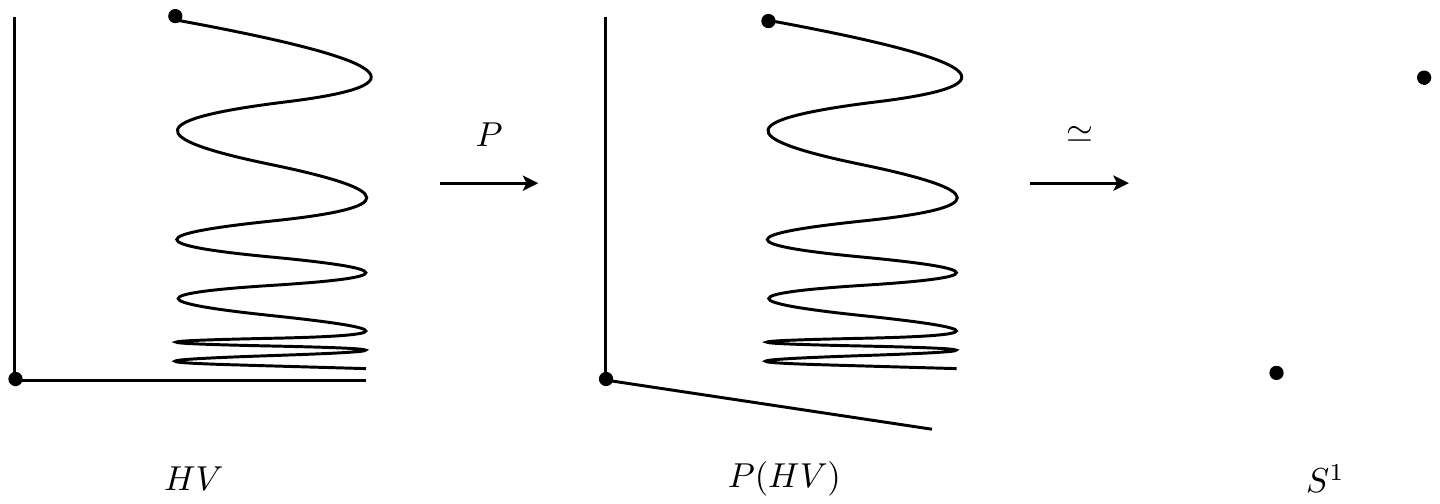}
    \caption{Schematic representation of Peanification and strong deformation retract. Presented are the intersections with $\f\in \{0,\pi\}$.}
\end{figure}

\subsection{Braiding  $HV$'s}

The next step is to take a countable union of Harmonic Vases, each of which will correspond to one generator of the group $G$. While doing so we have to be careful to maintain compactness, inner-heights of every $HV$ and to avoid intersections of different vases except at pedestal and at $\f=0$. Compactness will be preserved by decreasing the height of vases (namely decreasing $m$). Inner-heights will be different for suitable choice of parameters $p$ (namely they have to be algebraically independent over $\QQ$). Intersections will  be avoided using additional Euclidean dimension. Let us first introduce the notation that will explain how $HV's$ are embedded in $\RR^4$.

\begin{dfn}
\label{HV}
 \textbf{Harmonic Vase} with parameters $m, p \in \RR^+$ and $w\colon [-\pi,\pi]\times(0,m]\to \RR$ [notation: $HV(m,p,w)$] is the subset of $\RR^4$ defined as the union of two sets:
\begin{itemize}
  \item the pedestal $B(3,0)\cap (\RR^2 \times \{0\}^2)=\{(x,y,0,0)\in \RR^4, x^2+y^2\leq 9\},$ and
  \item the wall $W(m,p,w)$, parameterized as
  $$
  z\in (0,m],\quad \f\in [-\pi,\pi], \quad r:= \frac{|\f|}{\pi}\sin \frac{\pi p}{z}+2, \quad  w:=w(\f,z),
  $$
  where $(r,\f)$ are polar coordinates in $\RR^2\times \{0\}^2\subset \RR^4$, $z$ is the  coordinate representing $\{0\}^2\times \RR\times \{0\}$ so that $(r,\f,z)$ are cylindric coordinates in $\RR^3\times \{0\}$ and $w$ is the fourth coordinate representing $\{0\}^3\times \RR$.
\end{itemize}
\end{dfn}

We  define Braided Harmonic Vase ($BHV$) inductively. Let $\{p_i\}_{i\in \ZZ^+}$ be a sequence of positive numbers that are pairwise algebraically independent over $\QQ$ meaning that $p_i$ and $p_j$ are algebraically independent over $\QQ$ for every choice of $i\neq j$. To handle the intersections  let us describe them first. For $j<i; j,i\in \ZZ^+$ define
$$
H_i^j:=\{x\mid W(\frac{1}{i},p_i)\cap W(\frac{1}{j},p_j)\cap(\RR^2\times \{x\})\neq 0\}.
$$
In other words, $H_i^j$ is set of all heights where $W(\frac{1}{i},p_i)$ and $W(\frac{1}{j},p_j)$ intersect. Note that each of these sets is discrete in $(0,\fr{1}{i}]$: algebraic independence guarantees that no inner-height of $W(\frac{1}{i},p_i)$ is in $H_i^j$, inner-heights converge to $0$  and there are only finitely many elements of $H_i^j$ between two any two inner-heights. Consequently the finite union $H_i:=\cup_{j<i}H_i^j$ is discrete.
Hence there exist functions
$$
w_i\colon (0,\fr{1}{i}] \to [0,\fr{1}{i}]
$$
with the following properties:

\begin{equation}\label{Wcoord1}
w_i(x)<x\quad \forall x;
\end{equation}

\begin{equation}\label{Wcoord1.5}
 \quad w_i(x)\neq w_j(x) \quad \forall j<i, \forall x\in H^j_i;
\end{equation}

\begin{equation}\label{Wcoord2}
\quad w_i \equiv 0 \emph{ on some neighborhood of inner-heights}.
\end{equation}

As we already mentioned, the first condition maintains compactness, the second one allows us to avoid intersections and the third one preserves neighborhoods of inner-curves in $\RR^3$.

\begin{dfn} \label{BHV}
Let $\{p_i\}_{i\in \ZZ^+}$ be a sequence of positive numbers that are pairwise algebraically independent over $\QQ$ and let $w_i\colon (0,\fr{1}{i}] \to [0,\fr{1}{i}]$ be a set of functions satisfying $(\ref{Wcoord1})$ and $(\ref{Wcoord2})$. \textbf{Braided Harmonic Vase} with parameters $\{p_i,w_i\}_i$ [notation: $BHV(\{p_i,w_i\}_i)$] is
$$
\bigcup_{i\in \ZZ^+} HV(\fr{1}{i},p_i,|\f|w_i).
$$
The \textbf{wall} of $BHV$ is union of the walls of Braided $HV'$s. The \textbf{pedestal} of $BHV$ is the pedestal of any (every)  Braided $HV'$s.
\end{dfn}

Note that  the function $|\f|w_i$ allows us to avoid intersections between $HV'$s except at $\f=0$. We now need to prove that every $BHV$ is compact and calculate its fundamental group.

\begin{figure}
\includegraphics[bb=100 570 350 780]{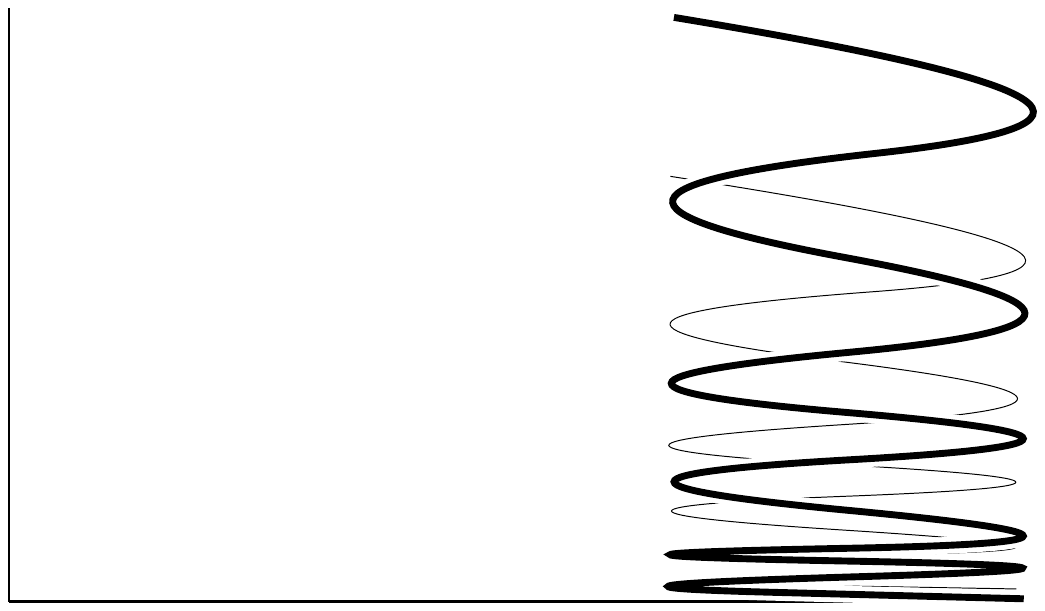}
\caption{Schematic representation of two incorporated $HV$'s in $BHV$ intersected with $\{\f=\pi\}$. Intersections are avoided using $w-$coordinate.}
\end{figure}

\begin{prop}\label{CompactBHV}
Every $BHV$ is compact.
\end{prop}

\proof
Take any Cauchy sequence $S$ in $BHV$. If $S$ converges to a point in $\RR^4$ with $z=0$ then because $w<z$ (the fourth coordinate is less that the third one by (\ref{Wcoord1})) that limit point is contained in  pedestal of $HV$. If $S$ converges to a point in $\RR^4$ with $z>0$ then we can assume that all points of $S$ have $z-$coordinate at least $\eps$ for some fixed $\eps>0$. Because $HV\cap \{z\geq \eps\}$ is compact (by the definition it is finite union of  images of $[\eps,\fr{1}{i}]\times [-\pi,\pi]$ as the heights of braided $HV'$s are decreasing) and contains $S$, it also contains the limit point.
\hfill $\blacksquare$
\bigskip

\remark Note that $BHV$ and $HV$ are  intersections of decreasing sequence of compacts, namely their closed $1/i-$neighborhoods. This fact gives alternative proof of compactness.

At this point we should emphasize the difference between a weak wedge and a metric wedge. Suppose $(X_i,x_i,d_i)_{i\in \ZZ^+}$ is a countable collection of pointed metric spaces with $x_i\in X_i$. Their \textbf{weak wedge} $\vee_{\ZZ^+}X_i$ is quotient space obtained by identifying all points $x_i$. Their \textbf{metric wedge} $\vee^m_{\ZZ^+}X_i$ is a metric  space obtained by identifying all points $x_i$ and defining the metric $d_\vee$ in the following way:
\begin{itemize}
  \item if $x,y\in X_j\subset \vee^m_{\ZZ^+}X_i$  for some $j\in \ZZ^+$ then $d_\vee(x,y):=d_j(x,y)$;
  \item else $d_\vee(x,y):=d_j(x,x_j)+d_k(y,x_k)$ where $x\in X_j\subset \vee^m_{\ZZ^+}X_i$, $y\in X_k\subset \vee^m_{\ZZ^+}X_i$.
\end{itemize}
The definition makes sense as $\vee^m_{\ZZ^+}X_i$ is pointwise union of sets $X_i.$ It is easy to prove that $d_\vee$ is indeed a metric.

In general $\vee_{\ZZ^+}X_i$ will almost never be metric because of the  topology at the wedge point. However, the topologies on natural subspaces $X_i$ are preserved by both wedges. Lemma \ref{wedge} proves that in many cases the homotopy types of maps from compact space to both wedges coincide.

\begin{dfn} \label{sdc}
Suppose $(X,x_0)$  is a pointed metric space. A \textbf{strong deformation contraction} of $X$  to $x_0$, is a continuous map $H\colon X\times I\to X$ so that
\begin{enumerate}
  \item $H|_{X\times \{0\}}=1|_X$;
  \item $H({X\times \{1\}})=H(\{x_0\}\times I)=x_0$;
  \item $d(H(x,t),H(y,t))\leq d(x,y), \forall t\geq 0$.
\end{enumerate}
\end{dfn}

\remark The use of this definition will be demonstrated in the proof of \ref{wedge} as  the metric wedge of strong deformation contractions is automatically continuous strong deformation contraction. Note that   metric wedge of strong deformation retractions needs not be a continuous map, it may fail to be  continuous at the wedge point.

\begin{lemma} Suppose that $R_i \colon X_i \to \{x_i\}$ are strong deformation contractions of pointed metric spaces $(X_i,x_i)$. Then the naturally defined map $R:=\vee_i R_i \colon \vee_i^m X_i \to \vee_i^m\{x_i\}$ on a metric wedge  is strong deformation contraction.
\end{lemma}

\proof We only need to show continuouity of $R$. Let $\{y_i\}_i=\{(p_i,t_i)\}_i$ be a Cauchy sequence of points in $\vee_i^m X_i\times I$ with limit $y=(p,t)$. If $y\in X_k-\{x_k\}\times I$ for some $k$ then $\lim_{i\to \infty}R(y_i)=R(y)$ as $R$ restricts to continuous $R_k$. If $y\in\vee_i^m\{x_i\}\times I$ then by definition \ref{sdc}
$$
d(D(y_i),D(y))=d(D(p_i,t_i),D(\vee_i^m\{x_i\},t_i))\leq d(p_i,\vee_i^m\{x_i\})\rightarrow_{i\to\infty}0.
$$
Hence $D$ is continuous.\hfill $\blacksquare$

\begin{lemma}\label{wedge} Let $r>0$ and suppose that for each $i\in \ZZ^+$ the $r-$neighborhood $U_i\subset X_i$ of point $x_i\in X_i$ in a metric space $(X_i,d_i)$ retracts to $x_i$ via strong deformation contraction $R_i$.
Then for each pointed compact space $(K,k_0)$ there is a natural bijection of homotopy classes of pointed maps $[K,\vee_{\ZZ^+}X_i]=[K,\vee^m_{\ZZ^+}X_i].$
\end{lemma}

\proof The natural map $\vee_{\ZZ^+}X_i\to \vee^m_{\ZZ^+}X_i$ is continuous hence we have a natural inclusion of the sets of maps $\mathcal{C}(K, \vee_{\ZZ^+}X_i)\subset \mathcal{C}(K,\vee^m_{\ZZ^+}X_i)$ and  $\mathcal{C}(K\times I, \vee_{\ZZ^+}X_i)\subset \mathcal{C}(K\times I,\vee^m_{\ZZ^+}X_i)$ which induce a well defined map $[K,\vee_{\ZZ^+}X_i]\to[K,\vee^m_{\ZZ^+}X_i]$. We will show that this map is a bijection.

First we prove that every map $f\colon(K,k_0) \to (\vee^m_{\ZZ^+}X_i,x_0)$ is homotopic $\rel k_0$ to a map $g\colon(K,k_0) \to (\vee^m_{S}X_i,x_0)\subset (\vee^m_{\ZZ^+}X_i,x_0)$ for some finite subset $S\subset \ZZ^+$. The finite metric and weak wedges coincide hence $g$ can naturally be considered as a map to $\vee_{\ZZ^+}X_i$.

Let $f\colon (K,k_0) \to (\vee^m_{\ZZ^+}X_i,x_0)$ be a map. The sets $X_i-U_i\subset \vee^m_{\ZZ^+}X_i$ are $2r$ disjoint hence there exists $n\in \NN$ so that compact $f(K)$ has empty intersection with $X_i-U_i$ for all $i\geq n$. Let $D$ be naturally defined homotopy on metric wedge
$$
X_1 \vee X_2 \vee \ldots \vee X_{n-1}\vee U_n \vee U_{n+1}\vee \ldots\subset \vee^m_{\ZZ^+}X_i
$$
so that $D(x,t):=x, \forall x\in X_1 \vee X_2 \vee \ldots \vee X_{n-1}$ and $D(x,t):=R_i(x,t)$ for all $x\in U_i, i\geq n$. Note that map $(x,t)\mapsto D(f(x),t)$ defined on $K\times I$ is a homotopy $\rel k_0$ between $f$ and the map $g$ whose image is contained in $\vee_{i< n}X_i$. Hence $g$ can be considered as a representative of $[f]$ in $[K,\vee_{\ZZ^+}X_i]$ which means that $[K,\vee_{\ZZ^+}X_i]\to[K,\vee^m_{\ZZ^+}X_i]$ is surjective. Using the same argument for space $K\times I$ we also see that $\vee_{\ZZ^+}X_i\to \vee^m_{\ZZ^+}X_i$ implies surjection on homotopies which means that $[K,\vee_{\ZZ^+}X_i]\to[K,\vee^m_{\ZZ^+}X_i]$ is injection hence bijection is proved.
\hfill $\blacksquare$

\begin{prop} \label{PeanifBHV}
$\pi_*(BHV)=\pi_*(\vee_{\ZZ^+}S^1).$
\end{prop}

\proof
Again the crucial step is to extract the space $P(BHV)$. The proof is almost the same as that of \ref{PeanifHV}

\begin{enumerate}
  \item Every point of  the  wall of $BHV$ has arbitrarily small simply connected neighborhood as  the wall itself is homeomorphic to a  countable union (in $\RR^4$) of $S^1\times (0,1]$ with common line $\{1\times I\}$. Hence  the topology of $P(BHV)$ at those points is no different from the topology in $BHV$.

  \item The point $(z=0,\f=0,r=2,w=0)$ (we will mark that point with $x_0$) also has arbitrarily small simply connected neighborhood. Let $\eps>0$ be very small and consider neighborhood $U_\eps$ of $x_0$ in $BHV$ that contains all points with $z<\eps$, $|\f|<\eps,$ $ 2+\eps >r>2-\eps.$ We will show that any point $x\in U_\eps$ can be connected to $x_0$ by a path which is enough for  the proof of our claim.

      If  the $z-$coordinate of $x$ equals $0$ then $x$ can be connected to $x_0$ because the pedestal (that contains both $x_0$ and $x$) is locally simply connected.

      If  the $z-$coordinate of $x$ equals $h\neq 0$ then $U_\eps \cap \{z=h\}$ is a finite wedge of open arcs, containing path from $x$ to a point $(z=zh,\f=0,r=2)$. This point can be connected to $x_0$ in $H_\eps$ by a straight line segment.

  \item The points on pedestal with $r>\frac{|\f|}{\pi}+2$ or $r<-\frac{|\f|}{\pi}+2$ are not limit points of the wall. Hence they all have arbitrarily small simply connected neighborhood.

  \item The points on pedestal with $\frac{|\f|}{\pi}+2 \geq r \geq-\frac{|\f|}{\pi}+2$ are limit points of the wall. But any point $x$, other that $x_0$ has  a neighborhood $U_x$  small enough, so that  the  path component of $U_x$ containing $x$ lies in  the  pedestal.
\end{enumerate}

Summing up, the only change of topology happens at the points in $(iv)$ which become separated from the wall: they have neighborhoods contained only in the pedestal. This means that $P(BHV)$  is homeomorphic to a wedge $B^2 \vee (\cup_{i\in \ZZ^+}W_i \cup \{x_0\})$ where $B^2$ represents  pedestal, wedge point represents $x_0$ and $(\cup_{i\in \ZZ^+}W_i)$ is the wall of $BHV$ (each $W_i$ represents the wall of some $HV$ braided in $BHV$).

Notice that the family $\{W_i\}_i$ is locally finite everywhere except at $x_0$. The union $\cup_i W_i$ in $P(BHV)$ can be replaced by a homeomorphic space: union of semi-open lateral sides of cylinders of increasing radius and decreasing height. To make the notation formal let $S(r,h)\subset \RR^3$ be semi-open lateral side of cylinder of radius $r$, height $h$ based at $(r,0,0)\in \RR^3:$
$$
S(r,h)=\{(x,y,z)\in \RR^3 \mid z\in (0,h]; d_{\RR^2}((x,y),(r,0))=r \}.
$$
Using this notation
$$
P(BHV)\cong \Big(B^2 \vee (\cup_{i\in \ZZ^+} S(2-\fr{1}{i},\frac{1}{i})\cup \{x_0\})\Big)
$$
where naturally $x_0=(0,0,0)\in \RR^3.$

Using  the obvious strong deformation retraction we see that
$$
P(BHV)\simeq V:=\cup_i\{(x,y)\in \RR^2 \mid d_{\RR^2}((x,y),(2-\fr{1}{i},0))=2-\fr{1}{i} \}.
$$
Using \ref{wedge} we get a natural bijection of homotopy classes of maps $[K,V]=[K,\vee_{\ZZ^+}S^1]$ for any compact space $K$. This bijection and \ref{MapsToPeano} imply $\pi_*(BHV)\cong \pi_*(\vee_{\ZZ^+}S^1 )$.
\hfill $\blacksquare$
\bigskip

\remark Note that  the space $V$ is not homeomorphic to a countable wedge of circles, it is homeomorphic to countable metric wedge of circles as the topology is not second countable.

\begin{figure}
\includegraphics[bb=130 550 370 770]{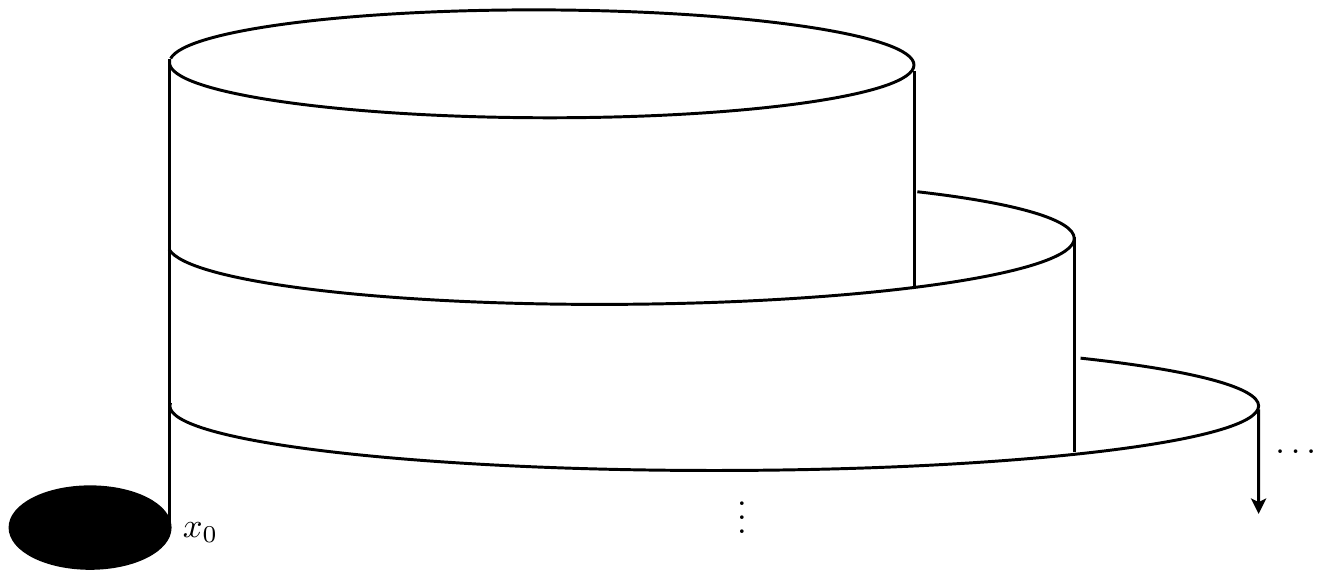}
\caption{$P(BHV)$.}
\end{figure}

\subsection{Attaching the relations}

We have constructed a compact metric space $BHV$ so that $\pi_1(BHV)=\langle g_1, g_2, \ldots \rangle$. In this subsection we will attach a disc $B^2$ to a space $BHV$ so that compactness will be preserved and the fundamental group will change to $\langle g_1, g_2, \ldots \mid r_1\rangle$. The following lemma explains how to attach a disc $B^2$ to $BHV$ at a certain "height" so that $BHV\cup B^2$ remains a subspace of $\RR^4$. Recall that $(z,w)$ stands for the pair of third and fourth Cartesian coordinates in $\RR^4$ respectively.
\bigskip

\notation For every $x,y\in \RR^+$ define:
\begin{itemize}
  \item $z_x:=(r=2,\f=0,z=x,w=0)$;
  \item $\g^x$ is  the linear path from $x_0=z_0$ to $z_x$;
  \item $\g^x_y$ is  the linear path from $z_y$ to $z_x$.
\end{itemize}
We will consider fundamental groups of $HV'$s and $BHV$ based at various points $z_h$. All the isomorphisms between differently based fundamental groups will be induced by paths $\g^*$ and $\g^*_*.$

\begin{lemma}\label{AttachRelation} Suppose $H= \cup_{i\in \ZZ^+} HV_i$ is a Braided Harmonic Vase where $HV_i$ are naturally incorporated Harmonic Vases. Let $h\in \RR^+$ and $r=[g_1 g_2 \ldots g_k]\in \pi_1(H,z_h)$ where each $[g_i]$ denotes one of two generators of some $[\g^h]^{-1}\ast\pi_1(HV_{j(i)},x_0)\ast [\g^h]$. Then  there exists $l\in \RR^+$ and an open topological $2-$disc $D$ so that:

\begin{enumerate}
  \item $D\subset \{h\geq z\geq l\}\subset \RR^4$
  \item $H\cap D=\emptyset$;
  \item $H \cup D \subset \RR^4$ is naturally homeomorphic to $ H\cup_r B^2$.
\end{enumerate}
\end{lemma}

\remark Parameters $h$ and $l$ allow us to attach $B^2$ to desired relation on $H$ as low (in terms of positive $z-$coordinate) as necessary. For this purpose the loop $g_1 g_2 \ldots g_k$ is based at $z_h$ as the loop along which we attach the disc should be contained in $\{h\geq z\geq l\}$.
\bigskip

\proof First we will define a path $\a$ in $H \cap  \{h\geq z\geq l\}$ so that:
\begin{itemize}
  \item $\a(0):=z_h$;
  \item $\a(1)=z_l$ for some $0<l\leq h$;
  \item $[\a\ast \g_l^h]=r\in \pi_1(H,z_h)$.
\end{itemize}
The construction of $\a$ is essentially a concatenation of two types of paths: vertical paths $\g^*_*$ (changing only  the $z-$coordinate) and generators of $\pi_1(W_i)$ near inner-heights.
\medskip

\textbf{Constructing the path}

Define $\a(0):=z_h$ and  let $a_1< h$ be an  inner-height of $HV_{j(1)}$. Define $\a_1(t):=\g_{z_{a_1}}^{z_h}(1-t)$. Note that  the image of $\a_1$ is contained in $H$. Appropriate orientation of a topological circle $HV_{j(1)}\cap \{z=a_1\}$ based at $z_{a_1}$ is a loop that  represents $[\g_{z_{a_1}}^h \ast g_1\ast (\g_{z_{a_1}}^h)^{-1}]\in \pi_1(H,z_{a_1})$. Let $\b'_1$ denote such loop based at $z_{a_1}$, i.e. $[\b'_1]=[\g_{z_{a_1}}^h \ast g_1\ast (\g_{z_{a_1}}^h)^{-1}]\in \pi_1(H,z_{a_1})$.

We still want to do a small correction of $\b'_1$. We want the function $t\mapsto \pi_z(\b'_1(t))$ to be decreasing where $\pi_z$ is projection to $z-$axis. The meaning of this condition will be explained later. Recall the meaning of  the function $w_{j(1)}$ from the definition of $BHV$ \ref{BHV}.  The  function $w_{j(1)}$ equals $0$ on a neighborhood $U_1$ of $a_1\in \RR$. It is not hard to see that we can homotope $\b'_1$ to another path (denote it by $\b_1$) in $HV_{j(1)}$ just by slightly changing $z-$coordinates within $U_1$ (decrease $z$ along $\b'_1$) so that we preserve starting point $\b_1'(0)=\b_1(0)$, make $t\mapsto \pi_z(\b_1(t))$ decreasing function and $[\b'_1]=[\b_1 \g_{\b_1(1)}^{a_1}]\in \pi_1(H,a_1)$.

We proceed by induction: let $a_2< \b_1(1)$ be an  inner-height of $HV_{j(2)}$. Define path $\a_2(t):=\g_{z_{a_2}}^{\b_1(1)}(1-t)$. The correct orientation of  the topological circle $HV_j(2)\cap \{z=a_2\}$ based at $z_{a_2}$ represents $[\g_{z_{a_2}}^h \ast g_2\ast (\g_{z_{a_2}}^h)^{-1}]\in \pi_1(H,z_{a_2})$. Let $\b'_2$ denote such loop based at $z_{a_2}$. Again we perturb $\b'_2$ to path $\b_2$ so that $t\mapsto \pi_z(\b_2(t))$ is decreasing and $\b'_2(0)=\b_2(0)$.

Having defined paths $\a_i$ (connecting paths) and $ \b_i$ (paths that represent $r_i$) for every $i\in \{1,\ldots,k\}$ we concatenate them  to get $\a$:
$$
\a:=\a_1 \ast \b_1 \ast \a_2 \ast \b_2 \ast \ldots \ast \a_k \ast \b_k
$$
Note that such defined $\a$ satisfies required conditions: $\a(0):=z_h$,  $\a(1)=z_l$ for some $0<l:=\b_k(1)$ and  $[\a\ast \g_l^h]=r\in \pi_1(H,z_h)$ by the construction. Also the map $t\mapsto \pi_z(\a(t))$ is decreasing.
\medskip

\textbf{Attaching the disc}

We will now attach a disc $B^2$ to $H$ along $\a$ (hence it's boundary will correspond to $r$)  so that the resulting space will still be embedded in $\RR^4$. First we define a map $f\colon \di I^2\to H$.

Define $f|_{\{0\}\times [1,1/2]}$ to be path $\a$ so that $f(0,0)=z_h$, define $f|_{\{0\}\times [1/2,1]}$ to be path $\g_{l}^{h}$ so that $f(0,1)=z_h$ and synchronize both parameterizations so that both paths are injective and $\pi_z f(0,1/2-t)=\pi_z f(0,1/2+t), \forall t\in [0,1/2]$. For every choice of $t\in [0,1]$ define $f(1,t):=(r=0,z=\pi_z(f(0,t)),w=0)$. Restating the definition,  left side of $\di I^2$ is path $\a \g_{z_l}^{z_h}$ and right side is the projection of $\a \g_{z_l}^{z_h}$ to axis $\{r=0,w=0\}$.

Define  the  map $f$ on  the  lower half of $I^2$ via straight line segments:
$$
f(s,t):=s(f(1,t))+(1-s)f(0,t)\qquad s\in (0,1), t\in [0,1/2].
$$

Note that $f(I\times \{s\})\subset \{z=\pi_z(f(0,s))\}$ for any $s\in [0,1/2]$. Because $t\mapsto \pi_z(\a(t))$ is decreasing this means that $f|_{I\times[1,1/2]}$ is injective.

We will use similar construction for the upper half of $I^2$ but we do want $f|_{(0,1)^2}$ to be injective. In order to ensure this we will use fourth dimension $w$. Let $g\colon I\times [1/2,1]\to I$ be a map with the following properties:
\begin{itemize}
  \item $g(\di(I \times [1/2,1]))=0$;
  \item $g(\Int(I \times [1/2,1]))> 0,$;
  \item $g(x,y)< \pi_z(f(1,y)), \forall x,y$.
\end{itemize}

The first condition will be required for the continuouity of $f$, the second allows $f$ to be injective on $(0,1)^2$, and the third one is necessary to maintain compactness of final construction.  Define  map $f$ on upper half of $I^2$ with perturbed straight line segments:
$$
f(s,t):=s(f(1,t))+(1-s)f(0,t)+ g(s,t)W \qquad s\in (0,1), t\in [1/2,1],
$$
where $W$ is unit vector along $w-$axis.

\begin{figure}
\includegraphics[bb=130 510 450 750]{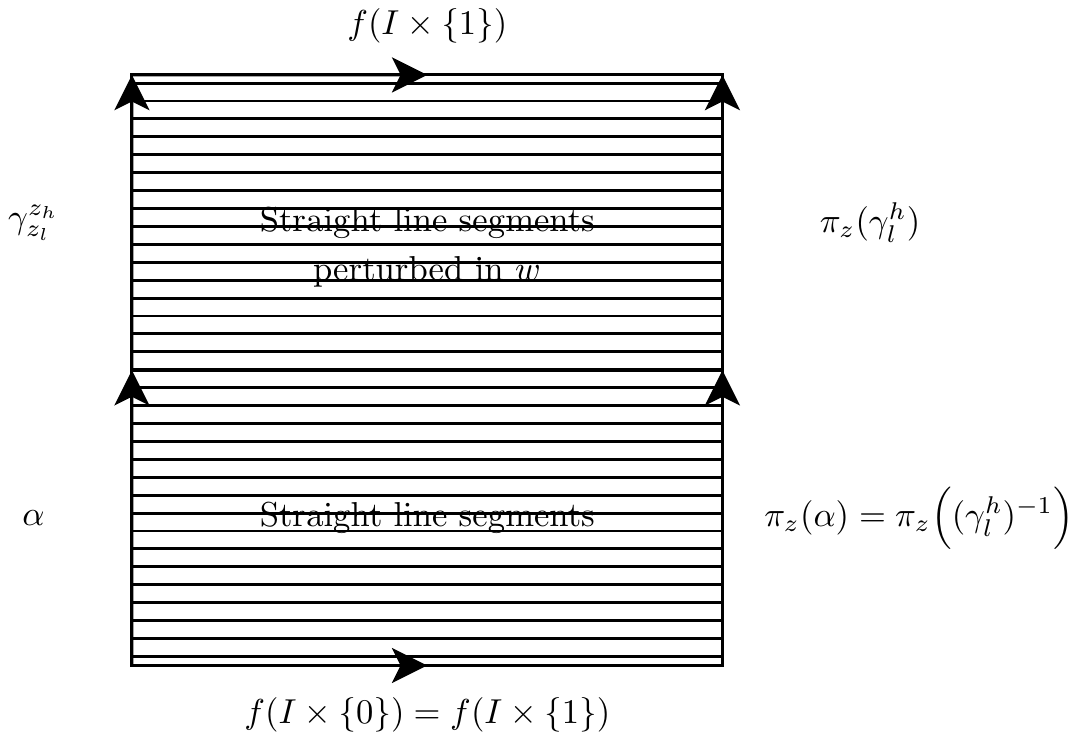}
\caption{Definition of map $f$ when attaching the disc.}
\end{figure}

Note  that $f|_{(0,1)^2}$ is injective. As mentioned above $f|_{I\times[1,1/2]}$ is injective due to the  map $t\mapsto \pi_z(\a(t))$ being decreasing. Also $f(I\times [0,1/2])\subset \{w=0\}=\RR^3\times\{0\}$ while all the points of $f((0,1)\times (1/2,1))$ have nontrivial $w-$coordinate by definition hence $f|_{(0,1)^2}$ is injective. Furthermore $f((0,1]\times [0,1])\cap H =\emptyset$. To see it let us analyze line segments $f(I\times \{t\})$. If  the $\f-$coordinate  $f((0,t))$ equals zero then $f((0,1]\times \{t\})\cap H =\emptyset$ which can be easily seen from figure \ref{HVcut1}. On other levels $f((0,t))$ is very close to an inner-height where by the definition $f((0,t))$ is the point with smallest $r-$coordinate in $H$ with given $(\f,z)-$coordinates. As the line segment ends in the point $(r=0,z=t,w=0)$ (so all the points of $f((0,1]\times \{t\})$ have even smaller $r-$coordinate) we see that $f((0,1]\times \{t\})\cap H = \emptyset$.

We now show that $f$ induces a map $f'$ on $B^2$ so that $[\di f'] = [r]$. Notice the equality of restrictions $f|_{I\times \{0\}} =f|_{I\times \{1\}}$. Also synchronized parametrization of $f|_{\{0\}\times I}$ implies  equality  $f(1,1/2-t)= f(1,1/2+t),$ $ \forall t\in [0,1/2]$. Identifying $(1,1/2-t)\sim (1,1/2+t), \forall t\in [0,1/2]$ and $(t,0)\sim (t,1), \forall t\in I$ we obtain quotient space $B^2$. It is easy to see that $f$ induces map $f'\colon B^2 \to \RR^4$ with the property $[\di f']=r \colon S^1\to H$. Also $f'|_{B^2-S^1}$ is injective and because $H$ is compact we obtain equality $H\cup f'(B^2)\cong H\cup_r B^2$.
\hfill $\blacksquare$
\bigskip

\remark Notice that $H$ is sufficiently nice to use Seifert Van-Kampfen theorem to obtain $\pi_1(H\cup_r B^2,x_0)=\langle g_1, g_2, \ldots \mid r\rangle$.

\subsection{Final construction}

Fix a countable group $G=\langle g_1, g_2, \ldots \mid r_1, r_2, \ldots \rangle$. Inductive use of \ref{AttachRelation} will provide us with a compact path connected subspace of $\RR^4$ that has $G$ as fundamental group.

\begin{thm} \label{Realization} For any countable group $G=\langle g_1, g_2, \ldots \mid r_1, r_2, \ldots \rangle$ there is a compact path connected subspace  $X_G\subset \RR^4$ so that $\pi_1(X_G,x_0)=G$.
\end{thm}

\proof We will define space $X_G$ inductively. Start with $X_0:=BHV$ and use \ref{AttachRelation} to attach $B^2$ to $X_0$ via map $r_1$ within $\{z\in (h_1,1) \}$  for some $h_1>0$ to get $X_1$. Proceed by induction: use \ref{AttachRelation} to attach $B^2$ to $X_k$ via map $r_{k+1}$ within $\{z\in (h_{k+1},h_{k}) \}$ for some $h_{k+1}>0$  to get $X_{k+1}$.  If there are only finitely many relations halt after finitely many steps, otherwise proceed with infinitely many steps and define $X_G:=\cup_i X_i$.

Space $X_G$ is natural subspace of $\RR^4$ and is path connected as every point of $X_i$ is path connected to $x_0\in X_i, \forall i.$ To prove $X_G$ is compact take any Cauchy sequence  $\{y_i\}_i$ in $X_G$. If $\pi_z(y_i)$ converge to $y_z> 0$ then use the fact that $X_G\cap \{z\geq y_z/2\}$ is compact as according to the construction only finitely many walls of braided $HV'$s in $X_0$ (heights of braided $HV's$ are decreasing) and only finitely many attached closed discs $B^2$ (see $(i)$ of \ref{AttachRelation}) intersect $\{z\geq y_z/2\}$ nontrivially. Both $HV'$s and discs $B^2$ are compact hence $X_G\cap \{z\geq y_z/2\}$ is compact therefore $\{y_i\}_i$ has a limit in $X_G$. If $\pi_z(y_i)$ converge to $0$ then note that $r-$coordinates of all elements are bounded by $3$ and $w-$coordinates of all elements are bounded by their $z-$ coordinates by the definition hence $w-$coordinate of limit point is $0$. Therefore the limit point of $\{y_i\}_i$ is contained in $\{r\leq 3,z=w=0\}$ which is the pedestal of $BHV$ hence contained in $X_G.$

The only thing left is to calculate $\pi_1(X_G,x_0)$. Again we will consider Peanification. Using the same argument as above  (see \ref{PeanifHV},\ref{PeanifBHV})  we see that Peanification of $X_G$ only moves pedestal apart from the wall of $H$ and relations, keeping it attached to the rest of the space only at $x_0$. Thus $P(X_G)$ is not compact.
We will prove that $[K,P(X_G)]=\cup_i[K,X_i], \forall K$ compact, where $[K,X_i]\subseteq [K,P(X_G)]$ is a subset of those homotopy classes of maps $K\to P(X_G)$ that have representative mapping $K\to X_i\subset P(X_G)$. This will mean $[K,X_G]=\cup_i[K,X_i], \forall K$ and hence (substituting $K$ for $S^1$ or $S^1\times I$) $\pi_1(X_G,x_0)=\langle g_1, g_2, \ldots \mid r_1, r_2, \ldots\rangle$ as every loop and every homotopy of $\pi_1(X_G,x_0)$ will be generated by some loop or homotopy of some $X_i$.
Notice that spaces $X_i$  are nice enough to use Seifert Van-Kampfen theorem and obtain $\pi_1(X_i,x_0)=\langle g_1, g_2, \ldots \mid r_1, r_2, \ldots,r_i\rangle$. Therefore the proof is concluded in the case of finitely many relations in the representation of $G$.

To prove equality $[K,P(X_G)]=\cup_i[K,X_i], \forall K$ for general countable group $G$ take any map $f\colon K \to P(X_G)$ and consider $P(X_G)\cap \{r=\f=w=0\}$.  For every $i\in \ZZ^+$ fix a point $x_i\in L_i:=B^2_i\cap \{r=\f=w=0\}\subset \{z\in (h_{i+1},h_{i})\}$, where $B^2_i$ is $B^2$ attached in $i^{th}$ step of construction of $X_G$. Recall from the definition that every  $B_i^2$ intersects $\{r=\f=w=0\}$.
Because $\lim_{k\to \infty}h_k=0$ the points $x_i$ are converging to  $\{z=r=\f=w=0\}\notin P(X_G)$ hence $f(K)$ can only hit finitely many  points $x_i$ which implies existence of $j\in \ZZ^+$ so that  $x_i\notin f(K), \forall i>j$.  Every point $x_i$ is contained in the interior of $B^2_i$ and because discs are apart from each other (separated by different zones of $z-$coordinate they occupy) there are natural strong deformation retractions $(B^2_i-\{x_i\})\to \di B^2_i\subset X_j, \forall i > j$ which induce homotopy of $f$ to a map $f'\colon K \to X_j$.
\hfill $\blacksquare$

\end{document}